\documentclass{amsart}
\usepackage{amsmath, enumerate, vmargin}

\newtheorem{thm}{Theorem}[section]
\newtheorem{prop}[thm]{Proposition}

\newtheorem{lem}[thm]{Lemma}
\newtheorem{defi}[thm]{Definition}
\newtheorem{remark}[thm]{Remark}
\newtheorem{example}[thm]{Example}
\newtheorem{pb}[thm]{Problem}

\newenvironment{rk}{\begin{remark}\rm}{\end{remark}}

\newcommand{\E}{{\mathcal E}}

\newcommand{\M}{{\mathcal M}}

\newcommand{\8}{\infty}
\newcommand{\el}{\ell}

\newcommand{\be}{\begin{eqnarray*}}
\newcommand{\ee}{\end{eqnarray*}}
\newcommand{\beq}{\begin{equation}}
\newcommand{\eeq}{\end{equation}}
\newcommand{\beqn}{\begin{equation*}}
\newcommand{\eeqn}{\end{equation*}}
\newcommand{\bsp}{\begin{split}}
\newcommand{\esp}{\end{split}}

\numberwithin{equation}{section}

\begin{document}
\title{{\bf Duality for symmetric Hardy spaces of noncommutative martingales}}
\author[T. N. Bekjan ]{Turdebek N. Bekjan}
\address{  L. N. Gumilyov Eurasian national University, Astana 010000, Kazakhstan}
\email{bekjant@yahoo.com}

\date{}
\begin{abstract}
We show the dual spaces of  conditional Hardy space and symmetric Hardy space of noncommutative martingales. We derive  relationship between  the symmetric Hardy space of noncommutative martingales and
its conditioned version.
\end{abstract}



\thanks{This work is supported by project 3606/GF4 of Science Committee of Ministry of Education and
Science of Republic of Kazakhstan.}

\maketitle


\section{Introduction}
The theory of noncommutative martingale inequalities has been rapidly developed. Many of the
classical martingale inequalities have been transferred to the noncommutative setting. Here we mention only three  of them
directly related with the objective of this paper. The first one is the noncommutative Burkholder-Gundy inequality. In the fundamental work \cite{PX}, Pisier/Xu established the noncommutative Burkholder-Gundy inequalities of noncommutative martingales in noncommutative $L_p$-spaces. In \cite{DPPS}, the authors proved Burkholder-Gundy inequalities for martingale difference
sequences in certain noncommutative symmetric spaces. Some versions can be found in  \cite{BCLJ2011,BC2012,DR}.

The second one is the noncommutative Burkholder-Rosenthal inequality. Junge/Xu \cite{JX2003,JX2007}) obtained  the noncommutative analogue of the classical Burkholder/Rosenthal
inequalities on the conditioned (or little) square function in noncommutative $L_p$-spaces. In \cite{DPPS}, the authors gave a generalization of Rosenthal inequality to noncommutative
symmetric spaces. Dirksen \cite{D} proved that if symmetric space's Boyd index satisfy $2<p_E\le q_E<\8$, then the Burkholder-Rosenthal inequalities hold in corresponding noncommutative symmetric space. Randrianantoanina and Wu \cite{RW} was proved  the Burkholder-Rosenthal inequalities hold if symmetric space has the Fatou property and its Boyd index satisfy $1<p_E\le q_E<2$.

The third  result is that independently Junge/Mei\cite{JM}  and Perrin\cite{P} obtained the following relationship between Hardy space of noncommutative martingales and
its conditioned version:
\beq\label{eq:relationship}
H^c_{p}(\M)=h^c_{p}(\M)+h^d_{p}(\M)
\eeq
for all $1\le p\le 2$.

The main novelty of our paper is the following dualities  for Hardy spaces and its conditional version:
If  $E$ is a separable symmetric space with $1< p_E\le q_E<2$ and $(E^{\times(\frac{1}{2})})^\times$ is separable, then
$$
(H^c_{E}({\mathcal{M}}))^*=H^c_{E^\times}({\mathcal{M}}),\quad (H^r_{E}({\mathcal{M}}))^*=H^r_{E^\times}({\mathcal{M}})$$
and
$$(h^c_{E}({\mathcal{M}}))^*=h^c_{E^\times}({\mathcal{M}}),\quad (h^r_{E}({\mathcal{M}}))^*=h^r_{E^\times}({\mathcal{M}}).
$$
As  applications of this result, we extended  \eqref{eq:relationship} to the noncommutative symmetric case. We use this result and Burkholder-Gundy inequalities to  obtain that  the  Burkholder-Rosenthal inequalities:
\begin{equation}\label{Burkholder-Rosenthal inequality}
L_{E}({\mathcal{M}})=h^c_{E}({\mathcal{M}})+h^d_{E}({\mathcal{M}})+h^r_{E}({\mathcal{M}})
\end{equation}
hold.

The organization of the paper is as follows. In Section $1$, we give some preliminaries and notations
on noncommutative
martingales and the noncommutative Hardy spaces. We prove the main results in Section $2$.

\section{Preliminaries}

Now let $E$ be a  quasi-Banach lattice. Let $0<\alpha<\infty$. $E$
is said to be $\alpha$-convex (resp. $\alpha$-concave) if there exists a
constant $C>0$ such that for all finite sequence $(x_{n})$ in $E$
$$
\begin{array}{c}
  \|(\sum|x_{n}|^{\alpha})^{\frac{1}{\alpha}}\|_E\le C (\sum\|x_{n}\|_E^{\alpha})^{\frac{1}{\alpha}},\\
(\mbox{resp}.\;\|(\sum|x_{n}|^{\alpha})^{\frac{1}{\alpha}}\|_E\ge C^{-1}
(\sum\|x_{n}\|_E^{\alpha})^{\frac{1}{\alpha}}).
\end{array}
$$
The least such constant $C$ is called the $\alpha$-convexity  (resp.
$\alpha$-concavity) constant of $E$ and is denoted by $M^{(\alpha)}(E)$ (resp.
$M_{(\alpha)}(E)$). For   $0<p<\infty,\; E^{(p)}$ will denote the
quasi-Banach lattice defined by
$$
E^{(p)}=\{x:\;|x|^{p}\in E\},
$$
equipped with the quasi-norm
$$
\|x\|_{E^{(p)}}=\||x|^{p}\|_{E}^{\frac{1}{p}}.
$$
It is easy to check that if $E$ is $\alpha$-convex and $\beta$-concave
then  $E^{(p)}$ is $p\alpha$-convex and $p\beta$-concave with
$M^{(p\alpha)}(E^{(p)})\le M^{(\alpha)}(E)^{\frac{1}{p}}$ and
$M_{(p\beta)}(E^{(p)})\le M_{(\beta)}(E)^{\frac{1}{p}}$. Therefore, if $E$ is
$\alpha$-convex  then  $E^{(\frac{1}{\alpha})}$ is $1$-convex, so it
can be renormed as a Banach lattice (see \cite{DDS,LT2,Xu1}).

 Let $([0£¬1],\Sigma, m)$ be the Lesbegue measure space and $L_0[0,1]$ be the space of all classes of Lebesgue measurable real-valued functions defined on $[0,1]$.  Let
 $x\in L_0[0,1]$. Recall that the
distribution function of $x$ is defined by
$$
\lambda_s(x) = m(\{\omega \in [0,1] :\; |x(\omega)| > s\}),\; s >0
$$
 and its decreasing rearrangement by
$$
\mu_t(x) = inf\{s > 0 : \; \lambda_s(x)\le t\},\;t > 0.
$$
For $x, y \in L_0[0,1]$ we say $x$ is submajorized by $y$, and write $x\preceq y$, if
$$
\int_{0}^{t}\mu_s(x)ds\le\int_{0}^{t}\mu_s(y)ds,\quad \mbox{  for all}\quad t > 0.
 $$

 By a symmetric  quasi Banach space on $[0,1]$ we mean a   Banach lattice $E$
 of measurable functions on $[0,1]$ satisfying the following
 properties: (i) $E$ contains all simple functions; (ii) if $x\in E$
 and $y$ is a measurable function on $[0,1]$ such that $\mu(y)\le \mu(x)$, then $y\in E$ and $\|y\|_{E}\le\|x\|_{E}$.

A symmetric quasi Banach space $E$ on $[0,1]$  is said to have the Fatou property if for every net $(x_{i})_{i\in I}$ in $E$ satisfying
$0\le x_i\uparrow$   and $\sup_{i\in I}\|x_i\|_E < \infty$ the supremum $x =\sup_{i\in I} x_i$ exists in $E$ and $\|x_i\|_E\uparrow\|x\|_E$. We say that $E$ has order continuous norm if for every net $(x_{i})_{i\in I}$ in $E$ such that $x_i\downarrow0$ we have $\|x_i\|_E\downarrow0$.

A symmetric  Banach  space $E$ on $[0,1]$ is called fully symmetric if, in addition, for $x \in L_0([0,1])$ and $y \in E$
with $x\preceq y$ it follows that $x \in E$ and $\|x\|_E \le \|y\|_E$.

The K$\ddot{o}$the dual of a symmetric Banach space $E$ on $[0,1]$ is the symmetric Banach
space $E^{\times}$ given by
$$
E^{\times} =\left\{x\in  L_0([0,1]):\; sup\{\int_{0}^1
|x(t)y(t)| dt\; : \|x\|_E \le 1\}< \infty\right\};
$$
$$
\|y\|_{E^{\times}} = sup\{\int_{0}^1
|x(t)y(t)| dt\; : \|x\|_E \le 1\},\quad
y \in  E^{\times}.
$$
The space $E^{\times}$ is  fully symmetric and has the Fatou property. It is isometrically
isomorphic to a closed subspace of $E^*$ via the map
$$
y\rightarrow L_y,\quad  L_y(x) =\int_{0}^1
x(t)y(t)dt\quad (x \in E).
$$
A symmetric Banach space $E$ on $[0,1]$ has the Fatou property
if and only if $ E = E^{\times\times}$ isometrically. It has order continuous norm if and only if
it is separable, which is also equivalent to the statement $E^{*} = E^{\times}$. Moreover, a
symmetric Banach  space which is separable or has the Fatou property is
automatically fully symmetric (cf. \cite{BS,B, KPS, LT2}).

   If E is a symmetric (quasi-)Banach  space  on $[0,1]$ which is
$\beta$-concave for some $\beta<\8$ then E has order continuous (quasi-)norm (cf. \cite{AB,DPPS}).

Let  $E_{i}$ be a symmetric quasi  Banach  space on
$[0,1],\; i=1, 2$. We define the pointwise product space $E_1\odot E_2$ as
\begin{equation}\label{multiplicativityspace}
    E_{1}\odot E_{2}=\{x:\;x=x_{1}x_{2},\;x_{i}\in E_{i},\;i=1,2\}
\end{equation}
with a functional $\|x\|_{E_1\odot E_2}$ defined by
$$
\|x\|_{E_1\odot E_2}=\inf\{\|x_{1}\|_{E_1}\|x_{2}\|_{E_2}:\;x=x_{1}x_{2},,\;x_{i}\in E_{i},\;i=1,2\}.
$$

By Theorem 2 in \cite{KLM}, we know that if $E_{i}$ is a quasi symmetric Banach  space on
$[0,1],\; i=1, 2$. Then $E_{1}\odot E_{2}$ is a quasi symmetric Banach  space on
$[0,1]$.

Using H$\ddot{o}$lder type inequality to Banach lattice (see Proposition 1.d.2 in \cite{LT2} or (3.1) in \cite{Xu1}) we obtain the following:
If $E$ is a symmetric  Banach  space on
$[0,1]$, then
\begin{equation}\label{squareroot}
E^{(\frac{1}{2})}=E\odot E.
\end{equation}

We use the following result( see (iii) of Theorem 1  and  Corollary 2 in \cite{KLM}).
\begin{thm}\label{pooductspace} Let $E,\;F$ be  symmetric Banach  spaces on
$[0,1]$.
\begin{enumerate}[\rm (i)]
  \item  If $0<p<\8$, then $(E\odot F)^{(p)}=E^{(p)}\odot F^{(p)}$.
  \item  If $1<p<\8$, then $(E^{(p)})^\times=(E^\times)^{(p)}\odot L_{p'}[0,1]$,
  where $\frac{1}{p}+\frac{1}{p'}=1$.
\end{enumerate}
\end{thm}
We also use the following  Lozanovskii  factorization theorem (Theorem 6 in \cite{Lo}).

\begin{thm}\label{lozanovskii}
Let $E$ be  a symmetric Banach space on $[0,1]$, then
$$
L_1[0,1]=E\odot E^{\times}.
$$

\end{thm}

\begin{prop}\label{lem:multi} Let $E$ be  a  symmetric Banach  space on $[0,1]$. If $E$ is $\alpha$-convex and $\beta$-concave for some $\alpha>1$ and $1<\beta<2$ then
 $ E =F\odot E^{\times}$,
 where $F=(E^{\times(\frac{1}{2})})^\times$ is separable.
\end{prop}

\noindent{\bf Proof.} Since $E^*=E^\times$, by Proposition 1.d.4 in \cite{LT2}, it follows that $E^\times$ is $\beta'$-convex and $\alpha'$-concave, where $\frac{1}{\alpha}+ \frac{1}{\alpha'}=1$ and $\frac{1}{\beta}+ \frac{1}{\beta'}=1$. Hence $E^{\times(\frac{1}{2})}$ is $\frac{1}{2}\beta'$-convex and $\frac{1}{2}\alpha'$-concave. Since $\frac{1}{2}\beta'>1$,  $E^{\times(\frac{1}{2})}$ can be renormed, with this equivalent norm, $E^{\times(\frac{1}{2})}$ becomes a  symmetric Banach space. On the other hand, from $ \frac{1}{2}\alpha'<\8$ follows that $E^{\times(\frac{1}{2})}$ is separable.

 Let $F=(E^{\times(\frac{1}{2})})^{\times}$.  Since $E^{\times(\frac{1}{2})}$  is separable,  $F=(E^{\times(\frac{1}{2})})^*$.  Therefore, $F$ is $\gamma$-concave, where $\gamma$ is conjugate of $\frac{1}{2}\beta'$ (i.e. $\frac{1}{\gamma}+ \frac{1}{\frac{1}{2}\beta'}=1$).  It follows that $\gamma<\8$, so $F$ is separable.

By  (ii) of Theorem \ref{pooductspace}, we have
 $$
 E=(E^\times)^\times=\big([E^{\times(\frac{1}{2})}]^{(2)}\big)^\times
 =\big([E^{\times(\frac{1}{2})}]^\times\big)^{(2)}\odot L^2[0,1]=F^{(2)}\odot L^2[0,1].
 $$
 Using Theorem \ref{lozanovskii} and (i) of Theorem \ref{pooductspace} we obtain that
  $$
 E=F^{(2)}\odot L^2[0,1]=F^{(2)}\odot[F\odot E^{\times(\frac{1}{2})}]^{(2)}=F^{(2)}\odot[F^{(2)}\odot E^{\times}]
 =F\odot E^{\times}.
 $$
 \hfill$\Box$

\subsection{Noncommutative symmetric  Banach spaces}

We use standard notation and notions from theory of noncommutative
$L^{p}$-spaces. Our main references are \cite{PX} and \cite{FK} (see
also \cite{PX} for more  historical references). Throughout this
paper, we denote by $\mathcal{M}$ a finite von Neumann algebra on
the Hilbert space $\mathcal{H}$ with a faithful normal normalized finite trace
$\tau$. The closed densely defined linear operator $x$ in ${\mathcal H}$
with domain $ D(x)$ is said to be affiliated with $\mathcal{M}$ if
and only if $u^{*}xu=x$ for all unitary operators $u$ which belong
to the commutant $\mathcal{M}'$ of $\mathcal{M}$. If $x$ is
affiliated with $\mathcal{M},$ then $x$ is said to be $\tau
$-measurable if for every $\varepsilon>0$ there exists a projection
$e\in \mathcal{M}$ such that $ e(\mathcal{H})\sqsubseteq D(x)$ and
$\tau( e^{\perp})<\varepsilon$ (where for any projection $ e$ we let
$e^{\perp}=1- e$ ). The set of all $ \tau $-measurable operators
will be denoted by $L_0(\mathcal{M})$. The set
$L_0(\mathcal{M})$ is a $\ast$-algebra with sum and product
being the respective closure of the algebraic sum and product.

 Let $x$ be a positive measurable operator and $x=\int_{0}^{\infty}
s\,de_s(x)$ its spectral decomposition. Composed with the spectral
measure $(e_s(x))$, the trace $\tau$ induces a positive measure
$d\tau(e_s(x))$ on $\mathbb{R}^+$. Thus we are reduced to the commutative case by
regarding $x$ as a function. In view of the preceding discussion
on functions, this justifies the following definition. Recall that
$e^\perp_s(x)=1_{(s,\infty)}(x)$ is the spectral projection of $x$
corresponding to the interval $(s,\infty)$.

 Let $x\in L_0(\mathcal{M})$. Define
 $$ \lambda_s(x)=\tau\big(e^\perp_s(|x|)),\quad s>0
 \quad\mbox{and}\quad
 \mu_t(x)=\inf\{s>0\;:\; \lambda_s(x)\le t\},\quad t>0.$$
The function $s\mapsto \lambda_s(x)$ is called the  distribution
function  of $x$ and the $\mu_t(x)$ the   generalized singular
numbers  of $x$.
For more details on generalized singular value function of measurable operators we refer to
\cite{FK}.

Let  $E$ be a symmetric  Banach space on
$[0,1]$. We define
$$
L_{E}({\mathcal{M}})=\{x  \in L_0(\mathcal{M}) :\;\mu_{.}(x)\in
E\};
$$
$$
\|x\|_{L_{E}({\mathcal{M}})}=\|\mu_{.}(x)\|_{E}, \quad \quad x
\in L_{E}({\mathcal{M}}).
$$
Then $(L_{E}({\mathcal{M}}),\|.\|_{L_{E}({\mathcal{M}})})$ is a  Banach
space (cf.  \cite{DDP1,S, Xu1}).

We need the following result (Theorem 5.6 in \cite{DDP3}, p. 745).

\begin{thm}\label{kotheduality} Let $E$ be  a separable symmetric Banach  space on $[0,1]$, then
$L_{E}(\mathcal{M})^* = L_{E^{\times}}(\mathcal{M})$ isometrically, with associated duality bracket given by
$$
\langle x, y\rangle = \tau(xy)\quad (x \in L_E(\mathcal{M}),\; y \in L_{E^{\times}}(\mathcal{M})).
$$
\end{thm}

In what follows, unless otherwise specified, we always denote by $E$  a symmetric  Banach space on
$[0,1]$.

Let $a=(a_{n})_{n\geq 0}$ be a finite sequence in
$L_{E}({\mathcal{M}})$, define
$$
\|a\|_{L_{E}({\mathcal{M}},\el_{c}^{2})}=\|(\sum\limits_{n\geq 0    }
|a_{n}|^{2})^{1/2}\|_{L_{E}({\mathcal{M}})},
$$
$$
\|a\|_{L_{E}({\mathcal{M}},\el_{r}^{2})}=\|(\sum\limits_{n\geq 0    }
|a_{n}^{*}|^{2})^{1/2}\|_{L_{E}({\mathcal{M}})}.
$$
This gives two  noms on the family of all finite sequences in
$L_{E}({\mathcal{M}})$. To see this, denoting by ${\mathcal{B}}(\ell_{2}) $ the
algebra of all bounded operators on $\ell_{2}$ with  its usual trace
$tr$, let us consider the von Neumann algebra tensor product
${\mathcal{M}}\otimes{ \mathcal{B}}(\ell_{2})$ with the product trace $\tau
\otimes tr.\;\tau \otimes tr$ is a semi-finite normal faithful
trace, the associated noncommutative $L_{E}$ space is denoted by
$L_{E}({\mathcal{M}}\otimes {\mathcal{B}}(\ell_{2})).$ Now, any finite sequence
$a=(a_{n})_{n\geq 1}$ in $L_{E}({\mathcal{M}})$ can be regarded as an
element in  $L_{E}({\mathcal{M}}\otimes {\mathcal{B}}(\ell_{2}))$ via the
following map
$$a\longmapsto T(a)= \left(
\begin{matrix}
a_{0} & 0 & \ldots \\
a_{1} & 0 & \ldots \\
\vdots & \vdots & \ddots
\end {matrix}
\right ),$$ that is, the matrix of $T(a)$ has all vanishing entries
except those in the first column which are the ${a_{n}}$'s. Such a
matrix is called a column matrix, and the closure in
$L_{E}({\mathcal{M}}\otimes {\mathcal{B}}(\ell_{2}))$ of all column matrices
is called the column subspace of $L_{E}({\mathcal{M}}\otimes
{\mathcal{B}}(\ell_{2}))$. Then
$$
\|a\|_{L_{E}({\mathcal{M}},\el_{c}^{2})}=\||T(a)|\|_{L_{E}({\mathcal{M}}\otimes
{\mathcal{B}}(\ell_{2}))}= \|T(a)\|_{L_{E}({\mathcal{M}}\otimes
{\mathcal{B}}(\ell_{2}))}.
$$
Therefore $\|.\|_{L_{E}({\mathcal{M}},\el_{c}^{2})}$ defines a  norm on the
family of all finite sequences of $L_{E}({\mathcal{M}}).$ The
corresponding completion is a  Banach space, denoted by
$L_{E}({\mathcal{M}},\el_{c}^{2})$. It is clear that if $E$ has the Fatou property, then a sequence
$a=(a_{n})_{n\geq 1}$ in $L_{E}({\mathcal{M}})$ belongs to
$L_{E}({\mathcal{M}},\el_{c}^{2})$ iff
$$
\sup_{n\geq 1}\|(\sum\limits_{k= 1}^{n}|a_{k}|^{2})^{1/2}\|_{L_E(\M)}<\infty.
$$
If this is the case, $(\sum\limits_{k= 1
}^{\infty}|a_{k}|^{2})^{1/2}$ can be appropriately defined as an
element of $L_{E}({\mathcal{M}}).$ Similarly, we may show that $
\|.\|_{L_{E}({\mathcal{M}},\el_{r}^{2})}$ is a  norm on the family of all
finite sequence in $L_{E}({\mathcal{M}}).$ As above, it defines a  Banach
space $L_{E}({\mathcal{M}},\el_{r}^{2}),$ which now is isometric to the
row subspace of $L_{E}({\mathcal{M}}\otimes {\mathcal{B}}(\ell_{2}))$
consisting of matrices whose nonzero entries lie only in the first
row. Observe that the column and row subspaces of
$L_{E}({\mathcal{M}}\otimes {\mathcal{B}}(\ell_{2}))$ are 1-complemented
subspace (by the definition of $E$ and Theorem 3.4 in \cite{DDP2}).

We also need $L_E^d(\M)$, the space of all sequences $a=(a_n)_{n\ge1}$ in $L_E(\M)$ such that
\be
\|a\|_{L_{E}^d({\mathcal{M}})}=\|diag(a_n)\|_{L_{E}({\mathcal{M}}\otimes
{\mathcal{B}}(\ell_{2}))}<\8,
\ee
where
$diag(a_n)$  is
the matrix with the $a_n$ on the diagonal  and zeroes
elsewhere.

\subsection{Noncommutative martingales}

Let ${\mathcal{M}}$ be a finite von Neumann algebra with a normalized normal faithful trace $\tau.$ Let
$({\mathcal{M}}_{n})_{n\geq 1}$ be an increasing sequence of von Neumann subalgebras of ${\mathcal{M}}$ such that $\cup_{n\geq 1}
{\mathcal{M}}_{n}$ generates ${\mathcal{M}}$ (in the $w^{*}$-topology). $({\mathcal{M}}_{n})_{n\geq 1}$ is called a
filtration of ${\mathcal{M}}.$ The restriction of $\tau$ to ${\mathcal{M}}_{n}$ is still denoted by $\tau.$ Let
${\mathcal{E}}_{n}={\mathcal{E}}(\cdot|{\mathcal{M}}_{n})$ be the conditional expectation of ${\mathcal{M}}$ with respect to
${\mathcal{M}}_{n}.$ Then ${\mathcal{E}}_{n}$ is a norm 1 projection of $L_{E}({\mathcal{M}})$ onto $L_{E}({\mathcal{M}}_{n})$ (see
Theorem 3.4 in \cite{DDP2}) and ${\mathcal{E}}_{n}(x)\geq 0$ whenever $x\geq 0.$

A noncommutative $L_{E}$-martingale with respect to $({\mathcal{M}}_{n})_{n\geq 1}$ is a sequence $x=(x_{n})_{n\geq 1}$
such that $x_{n} \in L_{E}({\mathcal{M}}_{n})$ and
\be
{\mathcal{E}}_n(x_{n+1})=x_n
\ee
for any $n \ge 1.$ Let $\|x\|_{L_E(\M)}=\sup_{n\geq 1}\|x_{n}\|_{L_E(\M)}.$ If $\|x\|_{L_E(\M)} <\infty,$ then $x$ is said to be a bounded $L_{E}$-martingale.

Let $x$ be a noncommutative martingale. The martingale difference sequence of $x,$ denoted by $dx=(dx_{n})_{n\geq 1},$ is defined as
\be
dx_{1}=x_{1},\quad dx_{n}=x_{n}-x_{n-1},\quad n\geq 2.
\ee
Set
\be
S^c_n (x)= \Big ( \sum_{k= 1  }^{n}|dx_{k}|^{2} \Big )^{\frac{1}{2}} \quad \mbox{and}
\quad S^r_n (x)= \Big ( \sum_{k=1}^{n}|dx_{k}^{*}|^{2} \Big )^{\frac{1}{2}}.
\ee
By the preceding discussion, if $E$ has the Fatou property, then $dx$ belongs to $L_{E}({\mathcal{M}},\el_{c}^{2})$ (resp. $L_{E}({\mathcal{M}}, \el_{r}^{2}))$ if and only if $(S^c_n (x))_{n \geq 1}$
(resp. $(S^r_n (x))_{n\geq 1}$) is a bounded sequence in $L_{E}({\mathcal{M}});$ in this case,
\be
S^c (x)= \Big ( \sum_{k= 1}^{\infty} | d x_{k} |^{2} \Big )^{\frac{1}{2}} \quad \mbox{and} \quad S^r (x)= \Big ( \sum_{k= 1 }^{\infty} | d x_{k}^* |^{2} \Big )^{\frac{1}{2}}
\ee
are elements in $L_{E}({\mathcal{M}}).$ These are noncommutative analogues of the usual square functions in the commutative martingale theory. It should be pointed out that the two sequences $S^c_n (x)\: \mbox{and}\: S^r_n (x)$ may not be bounded in
$L_{E}({\mathcal{M}})$ at the same time.

We define $H^c_{E}(\M)$ (resp. $H^r_{E}(\M)$) to be the space of all $L_{E}$-martingales such that $dx \in L_{E}(\M, \el_{c}^{2})$ (resp. $dx \in L_{E}(\M, \el_{r}^{2})$ ), equipped with the norm
\be
\|x\|_{H^c_{E}(\M)}=\|dx\|_{L_{E}(\M, \el_{c}^{2})
} \quad \big ( \mbox{resp.} \;
\|x\|_{H^r_{E}(\M)}=\|dx\|_{L_{E}(\M, \el_{r}^{2}) } \big ).
\ee
$H^c_{E}(\M)$ and $H^r_{E}(\M)$ are Banach spaces. Note that if $x \in H^c_{E}(\M),$
\be
\|x\|_{H^c_{E}(\M)} = \sup_{n \geq 0}\|S^c_n (x)\|_{L_{E}(\M)} = \|S^c (x)\|_{L_{E}(\M)}.
\ee
Similar equalities hold for $H^r_{E}(\M).$

We now consider the conditioned versions of square functions and Hardy spaces developed
in \cite{JX2003}. For a finite noncommutative $L_{E}$-martingale  $x=(x_{n})_{n\geq 1}$ define (with $\E_{0}=\E_1$)
\be
\|x\|_{h^c_{E}({\mathcal{M}})}=\| \Big ( \sum_{k\ge1}\E_{k-1}(|dx_{k}|^{2}) \Big )^{\frac{1}{2}}\|_{L_{E}({\mathcal{M }})} \ee
and
\be
 \|x\|_{h^r_{E}({\mathcal{M}})}= \|\Big ( \sum_{k\ge1}\E_{k-1}(|dx_{k}^{*}|^{2}) \Big )^{\frac{1}{2}}\|_{L_{E}({\mathcal{M }})}.
\ee
Let $h^c_{E}({\mathcal{M}})$ and $h^r_{E}({\mathcal{M}})$ be the corresponding completions. Then $h^c_{E}({\mathcal{M}})$ and $h^r_{E}({\mathcal{M}})$ are Banach
spaces. We define the column and row conditioned square functions as follows. For any finite
martingale $x = (x_n)_{n\ge1}$ in $L_E(M)$, we set
\be
s^c (x)= \Big ( \sum_{k\ge1 }\E_{k-1}(|dx_{k}|^{2}) \Big )^{\frac{1}{2}} \quad \mbox{and}
\quad s^r (x)= \Big ( \sum_{k\ge1}\E_{k-1}(|dx_{k}^{*}|^{2}) \Big )^{\frac{1}{2}}.
\ee
Then
\be
\|x\|_{h^c_{E}({\mathcal{M}})} = \|s^c (x)\|_{L_{E}({\mathcal{M }})} \quad \mbox{and}
\quad \|x\|_{h^r_{E}({\mathcal{M}})} = \|s^r (x)\|_{L_{E}({\mathcal{M }})}.
\ee
Let $x=(x_{n})_{n\geq 0}$ be a finite $L_{E}$-martingale, we set
$$
s^d (x)=diag(|dx_n|)
$$
We note that
\be
\|s^d (x)\|_{L_{E}({\mathcal{M}}\otimes
{\mathcal{B}}(\ell_{2}))}=\|dx_n\|_{L_{E}^d({\mathcal{M}})}
\ee
Let $h^d_E(\M)$ be the subspace of $L_{E}^d({\mathcal{M}})$ consisting of all martingale difference sequences.

\subsection{The space $L_E (\mathcal{M}; \el^{\8})$}

Recall that $L_E(\M; \el^{\8})$ is defined as the space of all sequences $(x_n)_{n \ge 1}$ in $L_E (\M)$ for which there exist $a, b \in L_{E^{(2)}}(\M)$ and a bounded sequence $(y_n)_{n \ge 1}$ in $\M$ such that $x_n = a y_n b$ for all $n \ge 1.$ For such a sequence, set
\beq\label{eq:E-MaxNorm}
\left \| ( x_n )_{n \ge 1} \right \|_{L_E (\M, \el^{\8})} : = \inf \big \{ \| a \|_{L_{E^{(2)}}(\M)} \sup_n\| y_n \|_{\8} \| b \|_{L_{E^{(2)}}(\M)} \big \},
\eeq
where the infimum runs over all possible factorizations of $(x_n)_{n \ge 1}$ as above. This is a norm and $L_E (\M; \el^{\8})$ is a Banach space (see \cite{D}).
As in \cite{JX2007}, we usually write
\be
\big \| {\sup_n}^+ x_n \big \|_E = \| ( x_n )_{n \ge 1} \|_{L_E (\M, \el^{\8})}.
\ee
We warn the reader that this suggestive notation should be treated with care. It is used for
possibly nonpositive operators and
\be
\big \| {\sup_n}^+ x_n \big \|_E \neq \big \| {\sup_n}^+ | x_n | \big \|_E
\ee
in general. However it has an intuitive description in the positive case: A positive sequence
$(x_n)_{n \ge1}$ of $L_E (\M )$ belongs to $L_E (\M; \el^{\8})$ if and only if there exists a positive $a \in L_E (\M)$ such
that $x_n \le a$ for any $n \ge 1$ and in this case,
\beq\label{eq:PositiveSequence}
\big \| {\sup_n}^+ x_n \big \|_E\le \inf \| a \|_{L_E (\M)}\le 2\big \| {\sup_n}^+ x_n \big \|_E,
\eeq
where the infimum runs over all possible positives $a \in L_E (\M)$ as above.  Indeed, let $(x_n)_{n \ge1}\in L_E (\M; \el^{\8})$. Then for $\varepsilon>0$, there exist $a, b \in L_{E^{(2)}}(\M)^+$ and a bounded sequence $(y_n)_{n \ge 1}$ in $\M$ such that $x_n = a y_n b$ for all $n \ge 1$, $\| a \|_{L_{E^{(2)}}(\M)}=\| b \|_{L_{E^{(2)}}(\M)}$, $ \sup_n\| y_n \|_{\8} =1$ and $\| a \|_{L_{E^{(2)}}(\M)}\| b \|_{L_{E^{(2)}}(\M)}<\| {\sup_n}^+ x_n \|_E+\varepsilon$. Define
the operator $c=(a^2+b^2)\in L_{E}(\M)^+$. Then there are
contractions $u,v\in M$ such that $a=c^{\frac{1}{2}}u,\; b=vc^{\frac{1}{2}}$ (see Remark 2.3 in \cite{DJ}). Hence,
$x_n =c^{\frac{1}{2}}{\rm Re}(uy_nv)c^{\frac{1}{2}}$, for all  $ n\ge 1$. It follows that $x_n \le c$  for all  $ n\ge 1$ and $\|c\|_{L_E (\M)}<2\| {\sup_n}^+ x_n \|_E+2\varepsilon$. So, the second inequality of \eqref{eq:PositiveSequence} holds. Conversely, if $x_n \le a$ for some $a\in L^+_E(M)$, then $x^{\frac{1}{2}}_n = u_na^{\frac{1}{2}}$ for a contraction $u_n\in\M$, and so $x_n =a^{\frac{1}{2}}u^*_nu_na^{\frac{1}{2}}$.
Thus $\| {\sup_n}^+ x_n \|_E\le \| a \|_{L_E (\M)}$. Hence, the first inequality of \eqref{eq:PositiveSequence} holds.

We define $L_E(\M; \el^1)$ to be the space of all sequences  $x=(x_{n})$ in  $L_{E}(\M)$ which
can be decomposed as
$$ x_{n}=\sum_{k=1}^{\infty}u_{kn}v_{nk}\quad (n\ge1)$$
for two families $(u_{kn})_{k,n\geq 1}$ and
$(v_{nk})_{n,k\geq 1}$ in $L_{E^{(2)}}(\M)$  satisfying
$$ \sum_{k,n=1}^{\infty}u_{kn}u_{kn}^{*}\in L_{E}(\M) \quad\mbox{and}\quad
 \sum_{n,k=1}^{\infty}v_{nk}^{*}v_{nk}
\in L_{E}(\M),
$$
where the series converge in norm. $L_E(\M; \ell^1)$ is a Banach space   when equipped with the norm
\be \|x\|_{L_{E}(\M;\ell^{1})}=\inf
\{\|\sum_{k,n=1}^{\infty}u_{kn}u_{kn}^{*}\|_{L_{E}(\M)}^{1/2}\|\sum_{n,k=1}^{\infty}v_{nk}^{*}v_{nk}\|_{L_{E}(\M)}^{1/2}\},
\ee where the infimum runs over all decompositions of $(x_n)$ as
above.
We will use the following fact (see \cite{D}). Let $x=(x_k)_{k\ge1}\in L_E(\M; \ell^1)$ for which $x_k\ge0$ for all
$k$. Then
\beq\label{normequalityofl1}
\|x\|_{L_{E}(\M;\ell^{1})}=\|\sum_{k\ge1} x_k\|_{L_E(\M)}.
\eeq

We need the following result (Theorem 5.3  and in Remark 5.4 in \cite{D}).

\begin{thm}\label{kothedualitymax} Let $E$ be  a separable symmetric Banach  space on $[0,1]$.
\begin{enumerate}[\rm(i)]
  \item If positive sequence
$(x_n)_{n \ge1}$ of $L_{E^{\times}} (\M )$ belongs to $L_{E^{\times}} (\M; \el^{\8})$,  then
\beq\label{eq:PositiveSequencenorm}
\big \| {\sup_n}^+ x_n \big \|_{E^{\times}} = \sup \big \{\sum_{k\ge1}\tau(x_ky_k) :\; y_k\in L_E (\M)^+,\;\|\sum_{k\ge1}y_k\|_{L_E (\M)}\le1 \big \}.
\eeq
  \item $L_{E}(\mathcal{M}; \el^1)^* = L_{E^{\times}}(\mathcal{M}; \el^{\8})$ isometrically, with respect to the duality bracket
$$
\langle x, y\rangle = \sum_{k\ge1}\tau(x_ky_k),
$$
where $x \in L_E(\mathcal{M}; \el^1)$ and  $y \in L_{E^{\times}}(\mathcal{M}; \el^{\8})$.
\end{enumerate}

\end{thm}


\section{Main results}

\begin{lem}\label{lem:opermulti} Let $E,E_1,E_2$ be  symmetric Banach  spaces on $[0,1]$ such that $ E =E_1\odot E_2$. If $x\in L_E(\M)^+$, then for $\varepsilon>0$, there exist $a\in L_{E_1}^+(\M)$ and $b\in L_{E_2}^+(\M)$ such that $x=ab,\; \|a\|_{L_{E_1}(\M)}\|b\|_{L_{E_2}(\M)}<\|x\|_{L_E(\M)}+\varepsilon$ and $a$ is invertible with bounded inverse.
\end{lem}

\noindent{\bf Proof.}
Let  $\mathcal{N}$ be the commutative von
 Neumann subalgebra of $\mathcal{M}$ generated by the spectral
 projection of $x$. Then $\mathcal{N}$ is isometric isomorphic
 to $L_\infty(\Omega,\Sigma,\mu)$ where $(\Omega,\Sigma,\mu)$
 a finite measure space. Hence, $x\in
 L_E(\mathcal{N})=L_E(\Omega,\mu)$. Since $ E =E_1\odot E_2$, for every $\varepsilon>0$, there are
 $x_1\in L_{E_1}(\Omega,\mu)^+=L_{E_1}(\mathcal{N})^+$ and
 $x_2\in L_{E_2}(\Omega,\mu)^+=L_{E_2}(\mathcal{N})^+$ such that $x=x_1x_2$ and
 $\|x\|_{L_{E}(\mathcal{N})}+\frac{\varepsilon}{2}>\|x_1\|_{L_{E_1}(\mathcal{N})}\|x_2\|_{L_{E_2}(\mathcal{N})}$. Let $\delta>0$. Set $a=x_1+\delta$ and $b=x_1(x_1+\delta)^{-1}x_2$.  Then $a\in L_{E_1}(\mathcal{N})\subset
L_{E_1}(\mathcal{M})$, $b\in L_{E_2}(\mathcal{N})\subset
L_{E_2}(\mathcal{M}),\;x=ab$ and $a$ is invertible with bounded inverse. For  enough small $\delta$ we have that
 $\|x\|_{L_{E}(\mathcal{M})}+\varepsilon>\|a\|_{L_{E_1}(\mathcal{M})}\|b\|_{L_{E_2}(\mathcal{M})}$.
Hence we obtain the desired result.
 \hfill$\Box$

\begin{thm}\label{duality}  Let $E$ be  a separable symmetric Banach  space on $[0,1]$ with $1<p_E\le q_E<2$.  If
$F=(E^{\times(\frac{1}{2})})^\times$ is separable, then we have
 \begin{enumerate}[\rm(i)]
 \item $(h^c_E(\M))^{*}=h_{E^\times}^c(\M) $ with equivalent norms.
 \item $(H^c_E(\M))^{*}=H_{E^\times}^c(\M)$ with equivalent norms.
 \end{enumerate}
 Similarly, $(h^r_E(\M))^{*}=h_{E^\times}^r(\M) $ and $(H^r_E(\M))^{*}=H_{E^\times}^r(\M) $ with equivalent norms.
\end{thm}
\noindent{\bf Proof.} (i) $1^\circ$ From $1<p_{E}\le q_{E}<2$, we obtain that $L_2(\M)\subset L_E(\M)$ with continuous inclusions, i.e. there exists a constant $C>0$ such that $\|x\|_E\le C\|x\|_{L_2(M)}$ for all $x\in L_2(\M)$. We  identify an element $x \in L_2(M)$ with the
martingale $(\E_n(x))_{n\ge1}$. By the trace-preserving property of conditional expectations
and the orthogonality in $L_2(\M)$ of martingale difference sequences, we get\\
$$\begin{array}{rl}
\|x\|_{h^c_{E}(\M)}&=  \|( \sum_{k\ge1}\E_{k-1}(|dx_{k}|^{2})  )^{\frac{1}{2}}\|_{L_{E}(\M)} \\
&\le C \|( \sum_{k\ge1}\E_{k-1}(|dx_{k}|^{2}) )^{\frac{1}{2}}\|_{L_{2}(\M)}\\
&=C\|x\|_{L_2(\M)},
\end{array}
$$
i.e. this martingale is in $h^c_{E}(\M)$.

 Let $y\in h_{E^\times}^c(\M)$. Since  $2<p_{E^\times}\le q_{E^\times}<\8$, it follows that $E^\times\subset L_2([0,1])$ with continuous inclusions. Hence, $\|y\|_{L_2(\M)}\le C_1 \|y\|_{h^c_{E}(\M)}<\8 $, i.e. $y$ is an $L_2$-
martingale.  Define $\phi_y$ by $\phi_y(x)=\tau(y^*x),\;\forall x\in L_2(\M)$. We must show that $\phi_y$ induces a continuous linear functional on $h^c_E(\M)$.

From the proof of Proposition \ref{lem:multi}, it follows that  $ E =F\odot E^{\times}$. Using Lemma \ref{lem:opermulti}, we obtain that  for $\varepsilon>0$, there exist $a\in L_{F}^+(\M)$ and $b\in L_{E^{\times}}^+(\M)$ such that $s^c(x)=ab,\; \|a\|_{L_{F}(\M)}\|b\|_{L_{E^{\times}}(\M)}<\|s^c(x)\|_{L_E(\M)}+\varepsilon$ and $a$ is invertible with bounded inverse.
Then, by the Cauchy-Schwarz inequality and the tracial property of $\tau$, we have

\be
\begin{array}{rl}
|\phi_y(x)|&= |\sum_{n\ge1}\tau(dy_n^*dx_n)|=|\sum_{n\ge1}\tau(\E_{n-1}(a)^{\frac{1}{2}}dy_n^*dx_n\E_{n-1}(a)^{-\frac{1}{2}})|\\
&\le\Big[\sum_{n\ge1}\tau(\E_{n-1}(a)^{\frac{1}{2}}|dy_n|^2\E_{n-1}(a)^{\frac{1}{2}})\Big]^{\frac{1}{2}}\Big[
\sum_{n\ge1}\tau(\E_{n-1}(a)^{-\frac{1}{2}}|dx_n|^2\E_{n-1}(a)^{-\frac{1}{2}}]^{\frac{1}{2}}\\
&=\Big[\sum_{n\ge1}\tau(a\E_{n-1}(|dy_n|^2))\Big]^{\frac{1}{2}}\Big[
\sum_{n\ge1}\tau(\E_{n-1}(a)^{-1}|dx_n|^2)\Big]^{\frac{1}{2}}\\
&=\rm{I}\cdot\rm{II}.
\end{array}
\ee
By  Theorem \ref{kotheduality}, we have
\be\begin{array}{rl}
\rm{I}^2 &=\sum_{n\ge1}\tau(a\E_{n-1}(|dy_n|^2))=\tau(a\sum_{n\ge1}\E_{n-1}(|dy_n|^2))\\
&\le\|a\|_{L_{F}(\M)}\|\sum_{n\ge1}\E_{n-1}(|dy_n|^2)\|_{L_{E^{\times(\frac{1}{2})}}(\M)}\\
&=\|a\|_{L_{F}(\M)}\|(\sum_{n\ge1}\E_{n-1}(|dy_n|^2))^\frac{1}{2}\|_{L_{E^{\times}}(\M)}^2\\
&=\|a\|_{L_{F}(\M)}\|y\|_{h_{E^\times}^c(\M)}^2.
\end{array}
\ee
To estimate II, we set
$$
s_{c,n}(x)=\big(\sum_{k=1}^n\mathcal{E}_{k-1}(|dx_k|^2)\big)^{\frac{1}{2}}, \quad \forall n\ge1
$$
and  $s_{c,0}(x)=|dx_1|$. Then
$$
\sum_{n\ge1}\E_{n-1}(|dx_n|^2)=\sum_{n\ge0}(s_{c,n+1}(x)^2-s_{c,n}(x)^2)).
$$
Applying Corollary 2.3 of \cite{C}  we obtain that $\E_{n-1}(a^{-1})\ge \E_{n-1}(a)^{-1}$. Hence, by  Theorem \ref{kotheduality}, we have
\be\begin{array}{rl}
{\rm II}^2&=\sum_{n\ge1}\tau(\E_{n-1}(a)^{-1}|dx_n|^2)\le\sum_{n\ge1}\tau(\E_{n-1}(a^{-1})|dx_n|^2)\\
&=\sum_{n\ge1}\tau(a^{-1}\E_{n-1}(|dx_n|^2))=\sum_{n\ge0}\tau(a^{-1}(s_{c,n+1}(x)^2-s_{c,n}(x)^2))\\
&=\tau(a^{-1}s^c(x)^2)=\tau(bs^c(x))\\
&\leq \|b\|_{L_{E^{\times}}(\M)}\|s^c(x)\|_{L_{E}(\M)}=\|b\|_{L_{E^{\times}}(\M)} \|x\|_{h^c_E(\M)}.
\end{array}
\ee
Combining the precedent estimations, for any finite $L_2$-martingale $x$, we deduce that
$$
|\phi_y(x)| \le \|y\|_{h_{E^\times}^c(\M)}\|x\|_{h^c_E(\M)}
$$
Thus $\phi$ extends to an element of $h^c_E(\M)^*$  with norm at most
$ \|y\|_{h_{E^\times}^c(\M)}$.

$2^\circ$ Let $\phi\in h_E^c(\M)^*$ of
norm one. As $L_2(\M) \subset h_E^c(\M)$, it follows that $\phi$ induces
a continuous functional  $\tilde{\phi}$  on $L_2(M)$.  Consequently, $\tilde{\phi}$ is
given by an element $y$ of $L_2(\M)$,
$$
\tilde{\phi}(y)=\tau(y^*x),\quad \forall x\in L_2(\M).
$$
As finite $L_2$-
martingales are dense in $h^c_{E}(\M)$ and in $L_2(M)$, we deduce that $L_2(M)$ is dense in $h^c_{E}(\M)$. We have
 \beq\label{eq:norm}
 \|\phi\|_{h_E^c(\M)^*}=\sup_{x\in L_2(\M),\|x\|_{h_E^c(\M)}\le1}|\tau(y^*x)|\le1.
 \eeq
We want to show that  $y\in h_{E^\times}^c(\M)$ and $\|y\|_{h_{E^\times}^c(\M)}\le C$.

Set
 \be
 z_n=\E_{n-1}(|dy_{n}|^2),\quad \forall n\ge1.
 \ee
 Using Theorem \ref{kotheduality} we obtain that
 \be
 \begin{array}{rl}
   \|y\|_{h^c_{E^\times}(\M)}^2&=\|\sum_{n\ge1}\E_{n-1}(|dy_{n}|^2)\|_{L_{E^{\times(\frac{1}{2})}}(\M)}
   =\|\sum_{n\ge1}z_{n}\|_{L_{E^{\times(\frac{1}{2})}}(\M)}\\&=\sup\big\{\sum_{n\ge1}\tau(z_{n}a):\;a\in L_{F}^+(\M)\; \mbox{and}\; \|a\|_{L_{F}(\M)}\le1\big\} \\
    &=\sup\big\{\sum_{n\ge1}\tau(z_{n}\E_{n-1}(a)):\;a\in L_{F}^+(\M)\; \|a\|_{L_{F}(\M)}\le1\big\}.
 \end{array}
 \ee

 Let $a\in L_{F}^+(\M)$  and $\|a\|_{F}\le1$.  Let $b$ be the martingale defined as follows:
$$
db_n=dy_n\E_{n-1}(a),\quad \forall n\ge1.
$$
 Using \eqref{eq:norm} we obtain that
 $$
 \tau(y^*b)\le \|b\|_{h_E^c(\M)}.
 $$
We have that
\be
\tau(y^*b)=\sum_{n\ge1}\tau(|dy_n|^2\E_{n-1}(a))=\sum_{n\ge1}\tau(\E_{n-1}(|dy_n|^2)\E_{n-1}(a))=\sum_{n\ge1}\tau(z_n\E_{n-1}(a)).
\ee
On the other hand, by the definition of $b$, we find
 $$\begin{array}{rl}
  s^c(b)^2&=\sum_{n\ge1}\E_{n-1}(\E_{n-1}(a)|dy_n|^2\E_{n-1}(a))\\
  &=\sum_{n\ge1}\E_{n-1}(a)\E_{n-1}(|dy_n|^2)\E_{n-1}(a)\\
  &=\sum_{n\ge1}\E_{n-1}(a)z_n\E_{n-1}(a).
 \end{array}
 $$
Let $c \in L^+_F (\M)$ such that  $\E_{n-1}(a) \le c$ for all $n\ge1$. Then for each $n$, there exists   a contraction $u_n\in\M$ such that  $\E_{n-1}(a)^{\frac{1}{2}}= u_nc^{\frac{1}{2}}$, and so
$$
  s^c(b)^2=\sum_{n\ge1}c^{\frac{1}{2}}u^*_n\E_{n-1}(a)^{\frac{1}{2}}z_n\E_{n-1}(a)^{\frac{1}{2}}u_nc^{\frac{1}{2}}=
  c^{\frac{1}{2}}\big(\sum_{n\ge1}u^*_n\E_{n-1}(a)^{\frac{1}{2}}z_n\E_{n-1}(a)^{\frac{1}{2}}u_n\big)c^{\frac{1}{2}}.
 $$
Whence
 $$
 \mu_t(s^c(b)^2)\le \mu_{\frac{t}{3}}(c^{\frac{1}{2}})\mu_{\frac{t}{3}}(\sum_{n\ge1}u^*_n\E_{n-1}(a)^{\frac{1}{2}}z_n\E_{n-1}(a)^{\frac{1}{2}}u_n)
 \mu_{\frac{t}{3}}(c^{\frac{1}{2}}),\quad \forall t\in [0,t].
 $$
On the other hand, from the proof of Proposition \ref{lem:multi}, we know that  $ E=F^{(2)}\odot L^2[0,1]$. By \eqref{squareroot}, we get
 $ E^{(\frac{1}{2})}=E\odot E=F^{(2)}\odot L^2[0,1]\odot L^2[0,1]\odot F^{(2)}=F^{(2)}\odot L^1[0,1]\odot F^{(2)}$.
Hence,
  \be
 \begin{array}{rl}
   \|b\|_{h_E^c(\M)}^2&=\|s^c(b)\|_{L_E(\M)}^2=\|s^c(b)^2\|_{L_{E^{(\frac{1}{2})}}(\M)}=\|\mu_t(s^c(b)^2)\|_{E^{(\frac{1}{2})}}\\
   &\le \|\mu_{\frac{t}{3}}(c^{\frac{1}{2}})
   \mu_{\frac{t}{3}}(\sum_{n\ge1}u^*_n\E_{n-1}(a)^{\frac{1}{2}}z_n\E_{n-1}(a)^{\frac{1}{2}}u_n)\mu_{\frac{t}{3}}(c^{\frac{1}{2}})\|_{E^{(\frac{1}{2})}}\\
   &\le  \|\mu_{\frac{t}{3}}(c^{\frac{1}{2}})\|_{F^{(2)}}\|\mu_{\frac{t}{3}}(\sum_{n\ge1}u^*_n\E_{n-1}(a)^{\frac{1}{2}}z_n\E_{n-1}(a)^{\frac{1}{2}}u_n)\|_1
   \|\mu_{\frac{t}{3}}(c^{\frac{1}{2}})\|_{F^{(2)}}\\
   &\le K\|c\|_F^{\frac{1}{2}}\big(\sum_{n\ge1}\tau(u^*_n\E_{n-1}(a)^{\frac{1}{2}}z_n\E_{n-1}(a)^{\frac{1}{2}}u_n)\big)\|c\|_F^{\frac{1}{2}} \\
    &\le K\|c\|_F\big(\sum_{n\ge1}\tau(z_n\E_{n-1}(a))\big),
 \end{array}
 \ee
 where $K$ is a constant depending only on norms of the delation operator  on the spaces $F^{(2)},\; L^1[0,1]$.
Using \eqref{eq:PositiveSequence}  and Corollary 5.4  of \cite{D1} we deduce that
 $$
 \|b\|_{h_E^c(\M)}\le C\|a\|_F^{\frac{1}{2}}\big(\sum_{n\ge1}\tau(z_n\E_{n-1}(a))\big)^{\frac{1}{2}},
 $$
 where $C$ is a constant.
 Combining the preceding inequalities, we obtain that
 $$
 \sum_{n\ge1}\tau(z_n\E_{n-1}(a))\le C^2\|a\|_F\le C^2.
 $$
 It follows that
 $\|y\|_{h_E^c(\M)}\le C$.
Thus we have finished the proof of (i).

(ii) $1^\circ$ Let $y\in H_{E^\times}^c(\M)$ and
define $\phi_y$ by $\phi_y(x)=\tau(y^*x),\;\forall x\in L_2(\M)$. We must show that $\phi_y$ induces a continuous linear functional on $H^c_E(\M)$. Let $x$ be a finite $L_2$-martingale such that $\|x\|_{H^c_E(\M)} < \8$.
\be
\begin{array}{rl}
|\phi_y(x)|&= |\sum_{n\ge1}\tau(dy_n^*dx_n)|=|\tau\otimes tr\big(\left(
                                                                   \begin{array}{ccc}
                                                                     dy_1 & 0 & \cdots \\
                                                                     dy_2 & 0 & \cdots\\
                                                                     \vdots &\vdots & \ddots \\
                                                                   \end{array}
                                                                 \right)^*\left(
                                                                   \begin{array}{ccc}
                                                                     dx_1 & 0 & \cdots \\
                                                                     dx_2 & 0 & \cdots\\
                                                                     \vdots &\vdots & \ddots \\
                                                                   \end{array}
                                                                 \right)
\Big)|\\
&\le \|\left(
\begin{array}{ccc}
dy_1 & 0 & \cdots \\
 dy_2 & 0 & \cdots\\
 \vdots &\vdots & \ddots \\
  \end{array}
   \right)^*\|_{L_{E^\times}({\mathcal{M}}\otimes
{\mathcal{B}}(\ell_{2}))}\|\left(
\begin{array}{ccc}
dx_1 & 0 & \cdots \\
dx_2 & 0 & \cdots\\
\vdots &\vdots & \ddots \\
\end{array}
\right)\|_{L_{E}({\mathcal{M}}\otimes
{\mathcal{B}}(\ell_{2}))}\\
&=\|y\|_{H^c_{E^\times}(\M)}\|x\|_{H^c_{E}(\M)}.
\end{array}
\ee
Thus $\phi_y$ extends to an element of $H^c_E(\M)^*$  with norm at most
$\|y\|_{H^c_{E^\times}(\M)}$.

$2^\circ$  Let $\phi\in H_E^c(\M)^*$ be such that $\|\phi\|_{H_E^c(\M)^*}\le1$. There exists  $y\in L_2(\M)$ such that
$$
\phi(y)=\tau(y^*x),\quad \forall x\in L_2(\M).
$$
By the density of $L_2(\M)$ in $H_E^c(\M)$, we have
 \beq\label{eq:normH}
 \|\phi\|_{H_E^c(\M)^*}=\sup_{x\in L_2(\M),\|x\|_{H_E^c(\M)}\le1}|\tau(y^*x)|\le1.
 \eeq
We use the same method in $2^\circ$ of the proof of (i) to show that  $y\in H_{E^\times}^c(\M)$ and $\|y\|_{H_{E^\times}^c(\M)}\le C$. By Theorem \ref{kotheduality}, we find that
\be
 \begin{array}{rl}
   \|y\|_{H^c_{E^\times}(\M)}^2&=\|\sum_{n\ge1}|dy_{n}|^2\|_{L_{E^{\times(\frac{1}{2})}}(\M)}
   \\&=\sup\big\{\sum_{n\ge1}\tau(|dy_{n}|^2a):\;a\in L_{F}^+(\M)\; \mbox{and}\; \|a\|_{L_{F}(\M)}\le1\big\} \\
    &=\sup\big\{\sum_{n\ge1}\tau(|dy_{n}|^2\E_{n}(a)):\;a\in L_{F}^+(\M)\; \|a\|_{L_{F}(\M)}\le1\big\}.
 \end{array}
 \ee
Let $a\in L_{F}^+(\M)$  and $\|a\|_{F}\le1$.  Set
$$
db_n=dy_n\E_{n}(a)-\E_{n-1}(dy_n\E_{n}(a)),\quad \forall n\ge1.
$$
Then  $b$ is  a martingale. By \eqref{eq:normH}, it follows that
 $$
 \tau(y^*b)\le \|b\|_{H_E^c(\M)}.
 $$
Since $(dy_n)_{n\ge1}$ is a martingale difference sequence, we obtain that
\be
\begin{array}{rl}
\tau(y^*b)&=\sum_{n\ge1}\tau(|dy_n|^2\E_{n}(a))-\sum_{n\ge1}\tau(dy_n^*\E_{n-1}(dy_n\E_{n}(a)))\\
&=\sum_{n\ge1}\tau(|dy_n|^2\E_{n}(a))-\sum_{n\ge1}\tau(\E_{n-1}(dy_n^*)dy_n\E_{n}(a))\\
&=\sum_{n\ge1}\tau(|dy_n|^2\E_{n}(a)).
\end{array}
\ee
Using the triangular inequality in $L_{E}({\mathcal{M}},\el_{c}^{2})$, we get
 \be
 \begin{array}{rl}
 \|b\|_{H_{E}^c(\M)}&=\|(db_n)_{n\ge1}\|_{L_{E}({\mathcal{M}},\el_{c}^{2})}\\
 &\le\|(dy_n\E_{n}(a))_{n\ge1}\|_{L_{E}({\mathcal{M}},\el_{c}^{2})}
 +\|(\E_{n-1}(dy_n\E_{n}(a)))_{n\ge1}\|_{L_{E}({\mathcal{M}},\el_{c}^{2})}.
\end{array}
 \ee
On the other hand, by Lemma 2.2 of \cite{Be}, there is a constant $r_E$ such that
 \be
 \|(\E_{n-1}(dy_n\E_{n}(a)))_{n\ge1}\|_{L_{E}({\mathcal{M}},\el_{c}^{2})}\le r_E\|(dy_n\E_{n}(a))_{n\ge1}\|_{L_{E}({\mathcal{M}},\el_{c}^{2})}.
 \ee
 Hence,
 \be
 \|b\|_{H_{E}^c(\M)}\le(1+r_E)\|(dy_n\E_{n}(a))_{n\ge1}\|_{L_{E}({\mathcal{M}},\el_{c}^{2})}.
 \ee
Let $c \in L^+_F (\M)$ such that  $\E_{n}(a) \le c$ for all $n\ge1$. Then for each $n$, there exists   a contraction $u_n\in\M$ such that  $\E_{n}(a)^{\frac{1}{2}}= u_nc^{\frac{1}{2}}$. As before, by
 $ E^{(\frac{1}{2})}=F^{(2)}\odot L^1[0,1]\odot F^{(2)}$, we have that
   \be
 \begin{array}{rl}
   \|(dy_n\E_{n}(a))_{n\ge1}\|_{L_{E}({\mathcal{M}},\el_{c}^{2})}^2
   &=\|c^{\frac{1}{2}}\big(\sum_{n\ge1}u^*_n\E_{n}(a)^{\frac{1}{2}}
   |dy_n|^2\E_{n}(a)^{\frac{1}{2}}u_n\big)c^{\frac{1}{2}}\|_{L_{E^{(\frac{1}{2})}}(\M)}\\
   &\le K\|c\|_F^{\frac{1}{2}}\big(\sum_{n\ge1}\tau(u^*_n\E_{n}(a)^{\frac{1}{2}}|dy_n|^2\E_{n}(a)^{\frac{1}{2}}u_n)\big)\|c\|_F^{\frac{1}{2}} \\
    &\le K\|c\|_F\big(\sum_{n\ge1}\tau(|dy_n|^2\E_{n}(a))\big).
 \end{array}
 \ee
By \eqref{eq:PositiveSequence}  and Corollary 5.4  of \cite{D1} we deduce that
 $$
\|b\|_{H_{E}^c(\M)}\le C\|a\|_F^{\frac{1}{2}}\big(\sum_{n\ge1}\tau(|dy_n|^2\E_{n}(a))\big)^{\frac{1}{2}},
 $$
 where $C$ is a constant.
 Combining the preceding inequalities, we obtain that
 $$
\tau(\sum_{n\ge1}|dy_n|^2a)= \sum_{n\ge1}\tau(|dy_n|^2\E_{n}(a))\le C^2\|a\|_F\le C^2.
 $$
 It follows that
 $\|y\|_{H_E^c(\M)}\le C$.
Thus we have finished the proof of (ii).

Passing to adjoint, we obtain the identities $(h^r_E(\M))^{*}=h_{E^\times}^r(\M) $ and $(H^r_E(\M))^{*}=H_{E^\times}^r(\M) $.
\hfill$\Box$

\begin{lem}\label{lem:duality} Let $E$ be  a separable symmetric Banach  space on $[0,1]$. Then we have
 $(h^d_E(\M))^*=h^d_{E^\times}(\M)$ with equivalent norms.
\end{lem}

\noindent{\bf Proof.}
Recall that $h^d_{E}(\M)$ consists of martingale difference sequences in $L_E^d(\M)$.
So $h^d_{E}(\M)$ is $2$-complemented in $L_E^d(\M))$  via the projection
$$P : \left\{\begin{array}{ccc}
              L_E^d(\M)) & \longrightarrow & h^d_{E}(\M)\\
		(a_n)_{n\geq 1} & \longmapsto & (\E_n(a_n)-\E_{n-1}(a_n))_{n\geq 1}
             \end{array}\right.$$
By Theorem \ref{kotheduality}, we obtain the desired result.
\hfill$\Box$

\begin{prop}\label{hardyinter} Let $E$ be  a separable symmetric Banach  space on $[0,1]$ with $2<p_E\le q_E<\8$. Then we have
 \begin{enumerate}[\rm(i)]
 \item $H^c_E(\M))=h_E^c(\M)\cap h^d_E(\M) $ with equivalent norms.
 \item $H^r_E(\M)=H_E^r(\M)\cap h^d_E(\M)$ with equivalent norms.
 \end{enumerate}
\end{prop}

\noindent{\bf Proof.} (i) Let $y\in H^c_E(\M))$. From the proof of Theorem 6.2 in \cite{D} we have that
\be
\|y\|_{H_E^c(\M)}\le \|y\|_{h_{E}^c(\M)\cap h_{E}^d({\mathcal{M}})}.
\ee
On the other hand, from the proof of Theorem 6.2 in \cite{D}, it follows that
\be
\|dy_n\|_{L_{E}^d({\mathcal{M}})}\le C_1 \|\sum_{n\ge1}r_n\otimes|dx_n|\|_{L_{E}(L^\8\overline{\otimes}\M)}.
\ee
Hence, by Theorem 4.7 in \cite{DPPS}, we get $\|dy_n\|_{L_{E}^d({\mathcal{M}})}\le C_2 \|y\|_{H_{E}^c(\M)}$. Since $p_{E^{(\frac{1}{2})}}=\frac{1}{2}p_E>1$, by Corollary 4.13 in \cite{D},  we have $\|y\|_{h_{E}^c(\M)}\le C_3 \|y\|_{H_{E}^c(\M)}$.
Thus
\be
 \|y\|_{h_{E}^c(\M)\cap h_{E}^d({\mathcal{M}})}\le C' \|y\|_{H_{E}^c(\M)}.
\ee

(ii) Passing to adjoint, we obtain the desired result.
\hfill$\Box$

\begin{prop}\label{hardysum} Let $E$ be  a separable symmetric Banach  space on $[0,1]$ with $1<p_E\le q_E<2$. If
$F=(E^{\times(\frac{1}{2})})^\times$ is separable, then we have
 \begin{enumerate}[\rm(i)]
 \item $H^c_E(\M))=h_E^c(\M)+h^d_E(\M) $ with equivalent norms.
 \item $H^r_E(\M)=h_E^r(\M)+ h^d_E(\M)$ with equivalent norms.
 \end{enumerate}
\end{prop}

\noindent{\bf Proof.} (i) By Proposition \ref{hardyinter}, it follows that there exist  constants $C>0$ such that
\beq\label{eq:inc1}
\max\{\|x\|_{h_{E^\times}^c(\M)},\;\|dy_n\|_{L_{E^\times}^d({\mathcal{M}})}\}\le C \|x\|_{H_{E^\times}^c(\M)},\quad \quad \forall x\in H_{E^\times}^c(\M).
\eeq
Let $y\in h^c_E(\M)$. By \eqref{eq:inc1} and Theorem \ref{duality},  we deduce that
\be
\|y\|_{H_{E}^c(\M)}=\sup_{\|x\|_{H_{E^\times}^c(\M)}\le1}|\tau(x^*y)|\le \sup_{\|x\|_{h_{E^\times}^c(\M)}\le C}|\tau(x^*y)|=C\|y\|_{h_{E}^c(\M)}.
\ee
Similarly,
$$
\|z\|_{H_{E}^c(\M)}\le C \|z\|_{h_{E}^d(\M)},\quad \forall z\in h_{E}^d(\M).
$$
Hence
\be
 \|y\|_{H_{E}^c(\M)}\le C\inf \{\|w\|_{h_{E}^d({\mathcal{M}})}+ \|z\|_{h_{E}^c(\M)}\},
\ee
where the infimum runs over all decomposition $y = w + z$ with $w$ in $h^c_{E}(\M)$ and $z$ in $h^d_{E}({\mathcal{M}})$. So, $h_E^c(\M)+h^d_E(\M)\subset H^c_E(\M) $. Using Theorem \ref{duality} and Lemma \ref{lem:duality} we obtain that $(h_E^c(\M)+h^d_E(\M))^*=h_{E^\times}^c(\M)\cap h^d_{E^\times}(\M)$. Sincer $E$ is  a separable symmetric Banach  space on $[0,1]$,   $E^*=E^\times$ is separable. On the other hand, $F=(E^{\times(\frac{1}{2})})^\times$ is separable, it follows that $E^\times$ is separable.   Applying Proposition \ref{hardyinter} we deduce that $H^c_E(\M))=h_E^c(\M)+h^d_E(\M) $ with equivalent norms.

(ii) Passing to adjoint, we obtain $H^r_E(\M)=h_E^r(\M)+ h^d_E(\M)$ with equivalent norms.
\hfill$\Box$

By Proposition \ref{hardyinter}, \ref{hardysum}  and Proposition 4.18 in \cite{DPPS},  we obtain the following:

\begin{thm}\label{normversion} Let $E$ be  a separable symmetric Banach  space on $[0,1]$ with  $1<p_E\le q_E<2$. If
$F=(E^{\times(\frac{1}{2})})^\times$ is separable, then we have
 \be
 L_E(\M)=h^c_E(\M)+h_E^d(\M)+ h^r_E(\M).
 \ee
\end{thm}

\begin{rk} Using  Theorem \ref{duality} and Theorem 6.2 in \cite{D}  we can obtain the result of Theorem \ref{normversion}.
\end{rk}

\begin{rk} Should Theorem \ref{normversion} be true whenever $\M$ is semi-finite and  $E$ has
the Fatou property. This was proved   in \cite{RW} by Randrianantoanina and   Wu via a very different method.
\end{rk}


\end{document}